\documentstyle[12pt]{amsart}
\input psfig

\begin{document}

\def\II{\mathaccent'27I} \def\1{\char'401}
\def\TT{\mathaccent'27T} \def\1{\char'401}
\def\DD{\mathaccent'27D} \def\1{\char'401}
\def\CC{\mathaccent'27C} \def\1{\char'401}
\def\BB{\mathaccent'27B} \def\1{\char'401}
\def\KK{\mathaccent'27K} \def\1{\char'401}
\def\XX{\mathaccent'27X} \def\1{\char'401}
\def\YY{\mathaccent'27Y} \def\1{\char'401}

\title[Local connectivity of the Mandelbrot set at 
certain points]{Local connectivity of the Mandelbrot set at
certain infinitely renormalizable points}
\thanks{1991 Mathematics Subject Classification, 
Primary : 58F23; Secondary: 30C10}
\author{Yunping Jiang}
\thanks{Partially supported by NSF grant and PSC-CUNY award}
\address{Department of Mathematics, Queens College of CUNY, 65-30
Kissena Blvd,\hfil\break Flushing, NY 11367}

\maketitle

\begin{abstract}
We construct a subset of the Mandelbrot set which is dense 
on the boundary of the Mandelbrot set and which consists of only 
infinitely renormalizable points
such that the Mandelbrot set is locally connected at every point of this subset.
We prove the local connectivity by finding
bases of connected neighborhoods directly.
\end{abstract}

\vskip50pt
\centerline{\bf Contents}
\vskip10pt
\noindent 1. Introduction
\vskip5pt
\noindent 2. Construction of bases of connected neighborhoods

\vskip 48pt

\noindent {\bf 1. Introduction}

\vskip5pt
Consider a quadratic polynomial $P_{c}(z)=z^{2}+c$, where $c$ is
a complex parameter. The filled-in Julia set
$K_{c}$ of $P_{c}$ is the set of points $z$ such that the orbit
$O(z) =\{ P^{\circ n}_{c}(z)\}_{n=0}^{\infty}$ is bounded. The Julia set
$J_{c}$ of $P_{c}$ is the boundary of $K_{c}$. 
The Mandelbrot set ${\cal M}$ is the set of parameters $c$
such that the critical orbit $O(c)=\{ P^{\circ n}_{c}(0)\}_{n=0}^{\infty}$
is bounded (see Fig. 1). Equivalently, 
${\cal M}$ is the set of parameters $c$ such that $J_{c}$ is
connected. For a parameter $c$ in the complement of the Mandelbrot set
${\cal M}$, the Julia set $J_{c}$ is a Cantor set. Douady and Hubbard (see
\cite{dh1}) proved that the Mandelbrot set ${\cal M}$ is connected.

\begin{figure}
\centerline{\psfig{figure=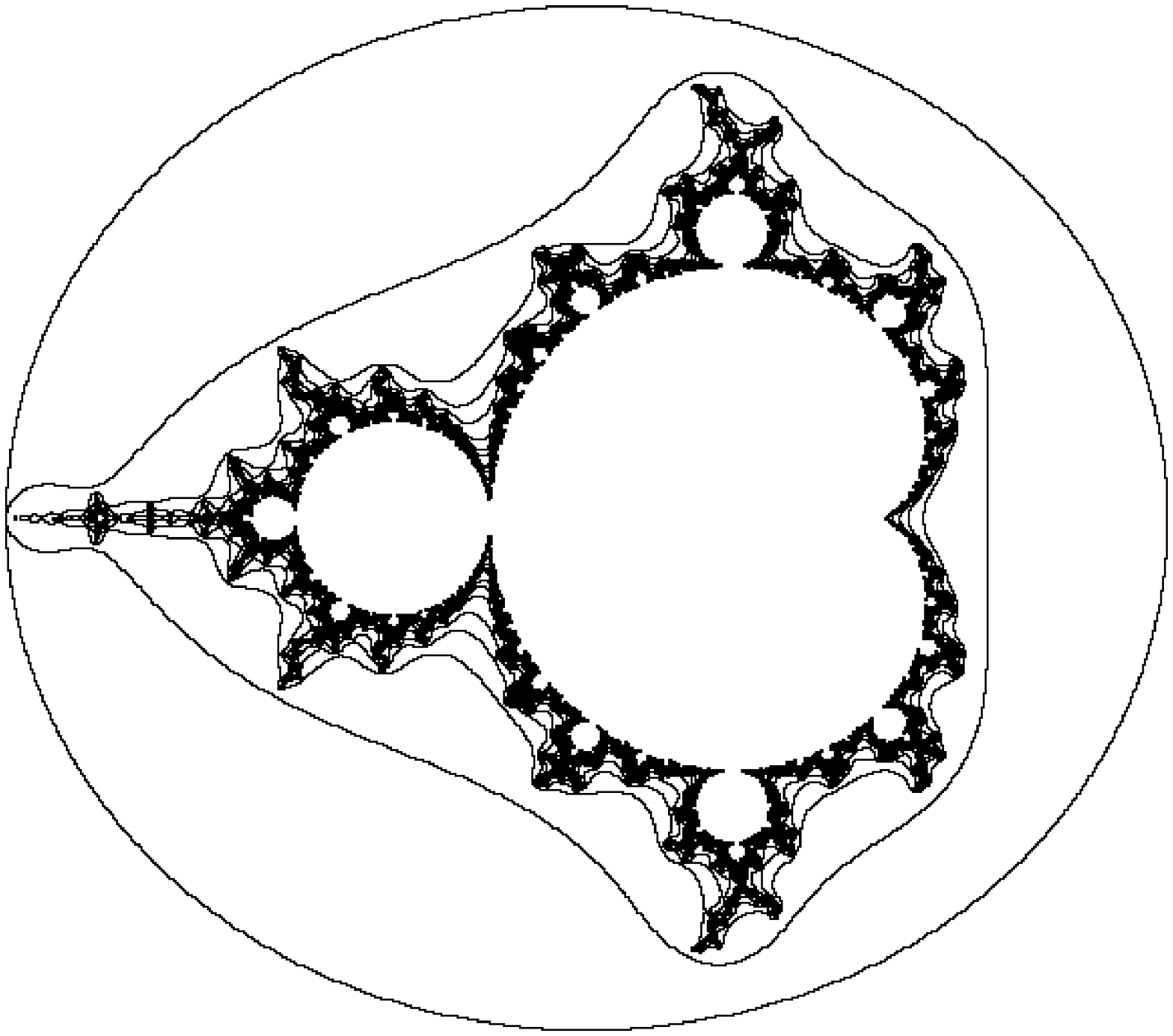,height=3in}}
\caption{Computer generating picture of the boundary of} 
\centerline{the Mandelbrot set and some equipotential curves}
\end{figure}

\vskip5pt
The Julia set of $P_{0}(z) =z^{2}$ is the round circle $S^{1}$. The map
$P_{0}$ restricted to $S^{1}$ is an expanding self-map. From the Structure
Stability Theorem in smooth dynamical systems (see \cite{pr,sh}), the
dynamical system $(P_{c}, J_{c})$ for any $|c|$ small enough is topologically
conjugate to $(P_{0}, S^{1})$. Thus $J_{c}$ for any $|c|$ small enough is a 
Jordan curve. It is also known that $J_{c}$ for any $|c|$ small enough is a
quasi-circle (see \cite{su}), and that the Hausdorff dimension $HD(J_{c})$
of $J_{c}$ is $1+|c|/(4\log 2) +\hbox{o}(|c|)$ (see \cite{ru}).
The techniques in the study of the dynamics of $P_{c}$ for any $|c|$ small
enough can be used to the study of the dynamics of a 
hyperbolic quadratic polynomial which is a quadratic polynomial
having an attractive or super-attractive periodic point in the complex plane. 
Let ${\cal HP}$ be the set of parameters $c$ such that $P_{c}$ is hyperbolic.
The main conjecture in this direction is that

\vskip5pt
\proclaim Conjecture 1. The set ${\cal HP}$ is open and dense in ${\cal M}$.

\vskip5pt
It is easy to verify that ${\cal HP}$ is open but to prove that ${\cal HP}$
is dense in ${\cal M}$ is a hard problem.

\vskip5pt
Douady and Hubbard \cite{dh2} proved that this conjecture follows from

\vskip5pt
\proclaim Conjecture 2. The Mandelbrot set ${\cal M}$ is locally connected.

\vskip5pt
In order to study the dynamics of quadratic polynomials more penetrating, 
Douady and Hubbard \cite{dh3} introduced the concept of a quadratic-like map
which is a proper, degree two, analytic branch cover $f: U\rightarrow V$ 
where $U$ and $V$ are open domains isomorphic
to a disc and $\overline{U}\subset V$. A quadratic polynomial 
$P_{c}$ is a quadratic-like map when restricted to
any domain bounded by an equipotential curve (see \cite{mi1}).

For a quadratic-like map $f: U\rightarrow V$, let $0$ be the unique 
branch point of $f$. The filled-in Julia set $K_{f}$ of $f$ is defined as 
$$
K_{f} =\cap_{n=0}^{\infty} f^{-n}(U).
$$
Suppose $K_{f}$ is connected (otherwise it is a Cantor set). 
The map $f$ is said to be once $n_{1}$-renormalizable if there are
subdomains $0\in U_{1} \subset U$ and $V_{1}\subset V$ isomorphic to
a disc such that 
$$
f_{1} =f^{\circ n_{1}} : U_{1} \rightarrow V_{1}
$$
is a quadratic-like map with connected filled-in Julia set
$K_{f_{1}}$. 
Otherwise, $f$ is called non-renormalizable.
If $f_{1}: U_{1}\rightarrow V_{1}$ is also $n_{2}$-renormalizable,
then we call $f$ is twice $(n_{1}, n_{2})$-renormalizable.
So on one can define a $k$-{\sl times} $(n_{1}, n_{2}, \cdots,
n_{k})$-renormalizable quadratic-like map as well as an infinitely
renormalizable quadratic-like map. 

Suppose $f: U\rightarrow V$ is a quadratic-like map. 
A point $p$ is called a periodic point of period $n\geq 1$ of $f$ 
if $f^{\circ n}(p)=p$ and $f^{\circ i}(p) \neq p$ for $1\leq
i<n$. The number $\lambda_{p} = (f^{\circ n})'(p)$ is called the
multiplier of $f$ at $p$. A periodic point $p$ of period $n$ of $f$ is
called parabolic if $\lambda_{p} =\exp (2\pi i p/q)$ where $p$ and $q$
are integers; irrationally indifferent if $\lambda_{p}=\exp (2\pi \theta i)$ where
$\theta$ is irrational; repelling if $|\lambda_{p}| >1$; attractive if
$0< |\lambda_{p}|<1$; super-attractive if $\lambda_{p}=0$. For a quadratic
polynomial $P_{c}(z)=z^{2}+c$, let $PCO(c)=\{ P_{c}^{\circ
n}(0)\}_{n=1}^{\infty}$ be its post-critical orbit.
We classify the parameters $c$ in ${\cal M}\setminus {\cal HP}$ as  
$$
{\cal P}= \{ c\in {\cal M}\; |\; P_{c} \hbox{ has a parabolic periodic point }\},$$
$$
{\cal I}=\{ c\in {\cal M} \; |\: P_{c} \hbox{ has an irrationally
indifferent periodic point }\},
$$ 
$${\cal N}=\{ c\in {\cal M} \; |\;  c \hbox{ is not in } {\cal P}\cup {\cal
I} \hbox{ and } 0\not\in \overline{PCO(c)}\; \}.
$$
$${\cal FR} =\{ c\in {\cal M}\; |\; c \hbox{ is not in } {\cal P}\cup {\cal
I}\cup {\cal N} \hbox{ and } P_{c} \hbox{ is not infinitely renormalizable
}\},$$
and 
$$
{\cal IR} =\{ c \in {\cal M} \; |\; P_{c} \hbox{ is infinitely
renormalizable }\}.
$$ 
In particular, a parameter $c$ in ${\cal N}$ is called Misiurewicz if
$PCO(c)$ consists of only finitely many points.
It is known that the set of all Misiurewicz points is dense on the boundary
of the Mandelbrot set (see \cite{cg}). A point $c$ in ${\cal IR}$ 
is called an infinitely renormalizable. 

Recent work of Yoccoz (see \cite{hu}) shows that the Mandelbrot set ${\cal
M}$ is locally connected at all points which are not infinitely
renormalizable. (We would like to note that the local connectivity 
of the Mandelbrot set ${\cal M}$ at all points in ${\cal P}$ is 
assured by many experts \cite{hu,sd} but a proof is still needed).
There remains the many infinitely renormalizable points for which to
complete the proof of Conjecture 2. However, we prove the following theorem 
to further confirm Conjecture 2.

\vskip5pt
\proclaim Main Theorem. There is a subset $\Upsilon$ of the Mandelbrot set ${\cal
M}$ which is dense on the boundary $\partial {\cal M}$ of the Mandelbrot set ${\cal
M}$ and which consists of only infinitely renormalizable points such
that the Mandelbrot set ${\cal M}$ is locally connected at every point 
$c$ of $\Upsilon$.

\vskip5pt
The proof of Main Theorem is in the next section. The talks with
Professors Tan Lei, C. Petersen, and A. Douady are very helpful.  
The article of J. Hubbard \cite{hu} and the book of Carleson and Gamelin
\cite{cg} are helpful during my understanding of
the structure of the Mandelbrot set ${\cal M}$ and the work of Yoccoz. 
I would like to thank them very much. I also would like to thank 
Professor M. Shishikura for pointing out the work of 
Eckmann and Epstein \cite{ee}.  Research is partially supported 
by NSF grant \# DMS-9400974 and PSC-CUNY award \# 6-65348. 
Research at MSRI is supported in part by NSF grant \# DMS-9022140.  
     
\vskip48pt

\noindent {\bf 2. Construction of bases of connected neighborhoods}

\vskip5pt
We give the proof of Main Theorem now. Let us start the construction of 
the subset $\Upsilon$ in Main Theorem around the
simplest point $-2\in {\cal M}$. Consider the polynomial $P(z)
=z^{2}-2$ and its Julia set $J=[-2, 2]$. 
It has a non-separate fixed point $\beta=2$ and a separate
fixed point $\alpha$, this menas that $J\setminus \{ \beta \}$ is still
connected and $J\setminus \{ \alpha\}$ is disconnected. 
Let $\gamma$ be a fixed equipotential curve of
$P$. Two external rays land at $\alpha$. Let $\Gamma$ be
the union of these two external rays. Let $D$ be the closed domain
bounded by $\gamma$. Let $D_{0}$ be the
domain containing $0$ and bounded by $\gamma$ and $P^{-1}(\Gamma)$ (see
Fig. 2). The restriction $P|({\bf C}\setminus (-\infty,-2])$ has two
inverse branches
$$
G_{0}: {\bf C}\setminus (-\infty,-2] \rightarrow {\bf LH}=\{ z=x+yi \in {\bf
C} \; |\; x<0\}
$$
and
$$
G_{1}: {\bf C}\setminus (-\infty,-2] \rightarrow
{\bf RH}=\{ z=x+yi \in {\bf C} \; |\; x>0\}.
$$
Let $D_{n} =G_{1}^{\circ n}(D_{0})$ and let $B_{n} =G_{0}(D_{n-1})$
for $n\geq 1$ (see Fig. 2). Since $2$ is an expanding fixed point of $P$ and
$P(-2) =2$, the diameter $\hbox{diam}(B_{n})$ tends to zero exponentially
as $n$ goes to infinity.

\addtocounter{figure}3
\begin{figure}
\centerline{\psfig{figure=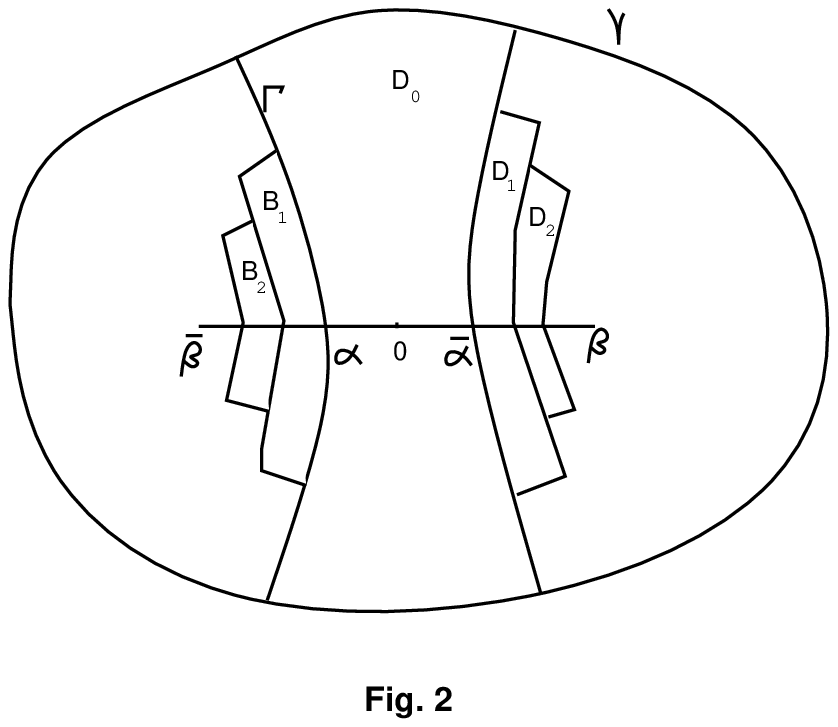}}
\end{figure}

Fig. 2 can be constructed as well for $c$ in a
small neighborhood about $-2$ because $\beta$,
$\alpha$, and $\Gamma$ are structurally stable.
Let $0\not\in U_{0}$ be a small neighborhood about $-2$ with
$\hbox{diam}(U_{0}) \leq 1$ such that Fig. 2 is
preserved for $c$ in $U_{0}$. Let $\beta (c)$, $\alpha (c)$,
$\Gamma (c)$, $D_{n}(c)$ for $0\leq n< \infty$, and $B_{n}(c)$ for
$1\leq n< \infty$ be the corresponding points, sets, and domains for
$P_{c}$, $c\in U_{0}$. Let $\overline{\beta}(c)\neq \beta (c)$ be
another pre-image of $\beta (c)$ under $P_{c}$. Then
the diameter $\hbox{diam}(B_{n})$ tends to $0$
and the set $B_{n}(c)$ approaches to $\overline{\beta}(c)$ as
$n$ goes to infinity uniformly on $U_{0}$.

\begin{figure}
\centerline{\psfig{figure=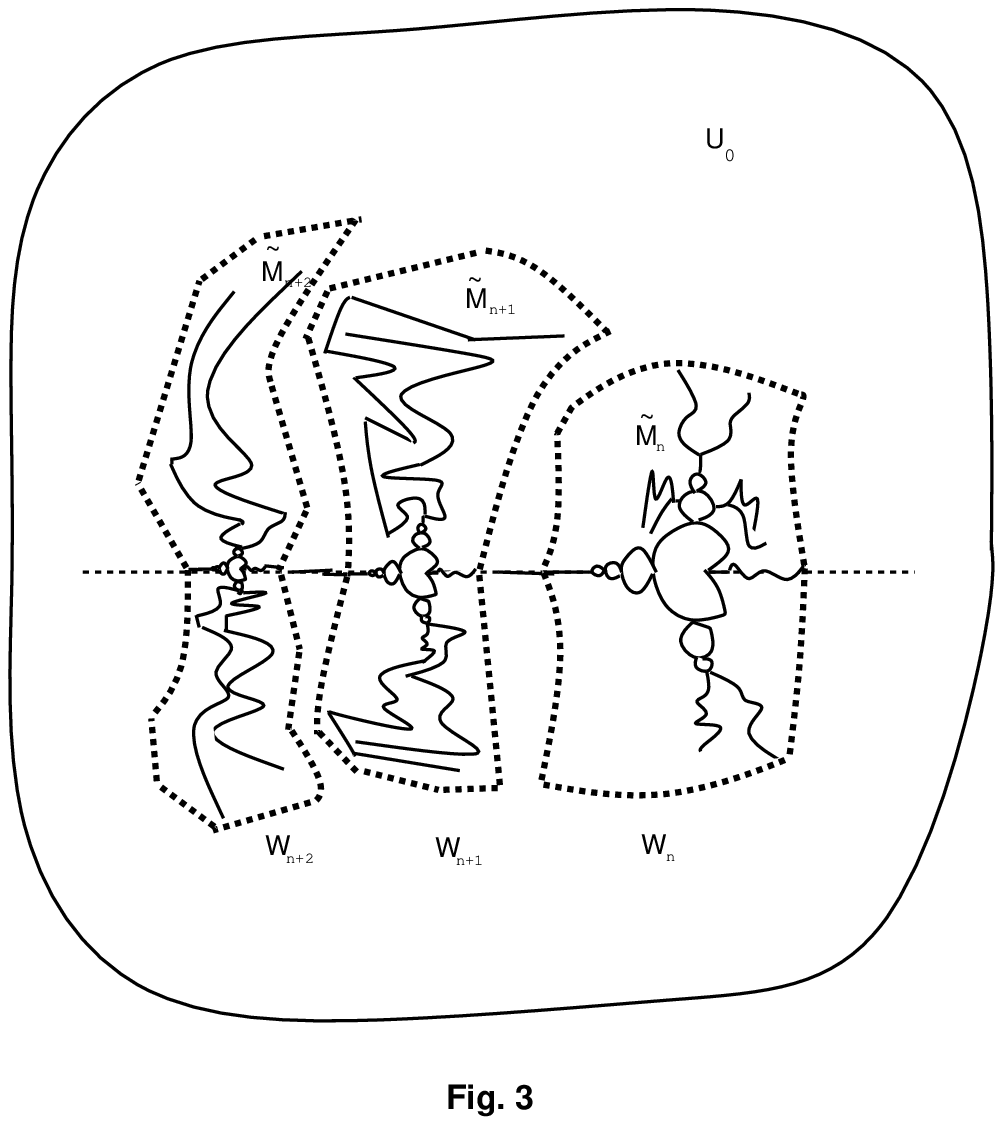}}
\end{figure}

Let
$$
W_{n} =\{ c \; |\; c\in B_{n}(c)\}.
$$
Because $W_{n}$ is bounded by external rays and an equipotential curve of
the Mandelbrot set ${\cal M}$, 
$$
\tilde{{\cal M}}_{n}=W_{n}\cap {\cal M}
$$ 
is connected.
Since the equation
$P_{c} (0)-\overline{\beta} (c)=0$
has a unique solution $-2$ in $U_{0}$,	the Rouch\'e Theorem assures that
$P_{c} (0)-x=0$ has a unique solution for any $x$ in $B_{n}(c)$ and large
$n$ and the solution is close to $-2$. Therefore, there is an integer
$N_{0}>0$ such that for $n\geq N_{0}$, 
$$
W_{n}\subset U_{0}.
$$
That is, 
$$
\hbox{diam}(W_{n})\leq 1
$$ 
for $n\geq N_{0}$ (see Fig. 3).

For each $W_{n}$ where $n\geq N_{0}$, since $c\in B_{n}(c)$,
$C_{n}(c)=P_{c}^{-1}(B_{n}(c))$ is a simply connected domain containing
$0$. It is also a subdomain of $D_{0}(c)$. Let
$$ 
F_{n,c}=P_{c}^{\circ (n+1)}: \CC_{n}(c) \rightarrow \DD_{0}(c)
$$
is a quadratic-like map (see Fig. 4).
The family 
$$
\{ F_{n, c}: \CC_{n}(c) \rightarrow \DD_{0}(c) \; |\; c\in W_{n} \}
$$ 
is full (see \cite{dh3}). So $W_{n}$ contains a copy ${\cal M}_{n}$
of the Mandelbrot set ${\cal M}$ (see \cite{dh3}). 
For $c\in {\cal M}_{n}$, 
the Julia set $J_{F_{n,c}}$ of $F_{n,c}: \CC_{n}(c) \rightarrow \DD_{0}(c)$ 
is connected. Therefore, for
$$c \in \Upsilon_{1} =\cup_{n\geq N_{0}}^{\infty} {\cal M}_{n},$$
$P_{c}$ is once renormalizable.

\begin{figure}
\centerline{\psfig{figure=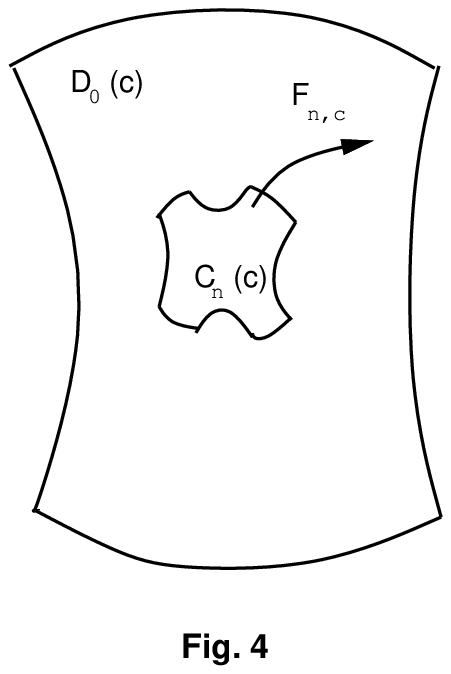}}
\end{figure}

For a fixed integer $i_{0}\geq N_{0}$, consider $W_{i_{0}}$ and
${\cal M}_{i_{0}}$; there is a parameter $c_{i_{0}} \in {\cal
M}_{i_{0}}$ such
that $$F_{i_{0}}=F_{i_{0}, c_{i_{0}}}: C_{i_{0}}=C_{i_{0}}(c_{i_{0}})
\rightarrow D_{i_{0}}=D_{0}(c_{i_{0}})$$
is hypbrid equivalent (see \cite{dh3}) to $P(z)=z^{2}-2$.
The quadratic-like map $F_{i_{0}}: C_{i_{0}}\rightarrow D_{i_{0}}$
has the non-separate fixed point
$\beta_{i_{0}}$ and the separate fixed point
$\alpha_{i_{0}}$. Let
$\overline{\beta}_{i_{0}}$ be another pre-image of $\beta_{i_{0}}$ under
$F_{i_{0}}$. Let $\Gamma_{i_{0}}$ be the external rays of
$P_{c_{i_{0}}}$ landing at $\alpha_{i_{0}}$. Let $D_{i_{0}0}$ be the
domain containing $0$ and bounded by $\partial C_{i_{0}}$ and
$F_{i_{0}}^{-1}(\Gamma_{i_{0}})$. Let $\overline{\beta}_{i_{0}} \in
E_{i_{0}0}$ and $\beta_{i_{0}} \in E_{i_{0}1}$ be the components
of the closure of $C_{i_{0}}\setminus D_{i_{0}0}$. Let $G_{i_{0}0}$ and
$G_{i_{0}1}$ be the inverses of $F_{i_{0}}|E_{i_{0}0}$ and
$F_{i_{0}}|E_{i_{0}1}$. Let
$$
D_{i_{0}n}=G_{i_{0}1}^{\circ n}(D_{i_{0}0})
$$
and let
$$ 
B_{i_{0}n} = G_{i_{0}0}(D_{i_{0}(n-1)})
$$
for $n\geq 1$.
Since $\beta_{i_{0}}$, $\alpha_{i_{0}}$, and $\Gamma_{i_{0}}$ are
structurally stable,
we can find a small neighborhood
$U_{i_{0}}$ about $c_{i_{0}}$ with $\hbox{diam}(U_{i_{0}})\leq 1/2$ such
that
the corresponding domains $B_{i_{0}}(c)$ can be constructed for $P_{c}$,
$c\in U_{i_{0}}$. Let
$$
W_{i_{0}n} =\{ c\in {\bf C} \; |\; F_{i_{0},c}(0) \in B_{i_{0}n}(c)
\}.
$$
The diameter $\hbox{diam}(B_{i_{0}n}(c))$ tends to zero
and the set $B_{i_{0}n}(c)$ approaches to
$\overline{\beta_{i_{0}}}(c)$
uniformly on $U_{i_{0}}$
as $n$ goes to infinity.
Since the equation $F_{i_{0},c}(0)-\overline{\beta}_{i_{0}}(c)=0$
has a unique solution $c_{i_{0}}$, the Rouch\'e Theorem implies that
$F_{i_{0}, c}(0)-x=0$ has a unique solution near $c_{0}$ for any $x\in
B_{i_{0}n}(c)$ and any $n$ large. Thus,
there is an integer $N_{i_{0}}\geq 0$ such that for $n\geq N_{i_{0}}$,
$W_{i_{0}n}\subseteq U_{i_{0}}$. Thus $\hbox{diam}(W_{i_{0}n})\leq 1/2$.
Since $W_{i_{0}n}$ for $n\geq N_{i_{0}}$ is bounded by external rays and
an equipotential curve of ${\cal M}$,
$\tilde{\cal M}_{i_{0}n}=W_{i_{0}}\cap
{\cal M}$ is connected. For each $c$ in $W_{i_{0}n}$, $n\geq N_{i_{0}}$,
let $C_{i_{0}n}(c) =F^{-1}_{i_{0},c}(B_{i_{0}n}(c))$.
Then
$$
F_{i_{0}n,c}=F^{\circ (n+1)}_{i_{0},c} : \CC_{i_{0}n}(c) \rightarrow
\DD_{i_{0}0}(c)
$$
is a quadratic-like map. Moreover,
$$ 
\{ F_{i_{0}n,c}\; |\; c\in W_{i_{0}n}\}
$$ 
is a full family. Thus
$W_{i_{0}n}$ contains a copy ${\cal M}_{i_{0}n}$ of the Mandelbrots set 
${\cal M}$. For
$$
c\in \Upsilon_{2} = \cup_{i_{0}\geq N_{0}} \cup_{i_{1}\geq
N_{i_{0}}} {\cal M}_{i_{0}i_{1}},
$$
$P_{c}$ is twice renormalizable.

We use the induction to complete the construction of the subset
around $-2$. Suppose we have constructed $W_{w}$ where
$w=i_{0}i_{1}\ldots i_{k-1}$ and $i_{0}\geq N_{0}$, $i_{1}\geq
N_{i_{1}}$, $\ldots$, $i_{k-1}\geq N_{i_{0}i_{1}\ldots i_{k-2}}$. Let
$v=i_{0}\ldots i_{k-2}$. There is a parameter $c_{w} \in {\cal M}_{w}$
such that $$F_{w}=F_{w, c_{w}}: C_{w}=C_{w}(c_{w})
\rightarrow D_{w}=D_{v0}(c_{w})$$
is hypbrid equivalent to $P(z)=z^{2}-2$.
The quadratic-like map $F_{w}: C_{w} \rightarrow D_{w}$ has the non-separate
fixed point $\beta_{w}$ and the separate fixed point $\alpha_{w}$. Let
$\overline{\beta}_{w}$ be another pre-image of $\beta_{w}$ under
$F_{w}$. Let $\Gamma_{w}$ be the external rays of
$P_{c_{w}}$ landing at $\alpha_{w}$. Let $D_{w0}$ be the
domain containing $0$ and bounded by $\partial C_{w}$ and
$F_{w}^{-1}(\Gamma_{w})$. Let $\overline{\beta}_{w} \in
E_{w0}$ and $\beta_{w} \in E_{w1}$ be the components
of the closure of $C_{w}\setminus D_{w0}$. Let $G_{w0}$ and
$G_{w1}$ be the inverses of $F_{w}|E_{w0}$ and
$F_{w}|E_{w1}$. Let
$$
D_{wn}=G_{w1}^{\circ n}(D_{w0})
$$
and let
$$ 
B_{wn} = G_{w0}(D_{w(n-1)})
$$
for $n\geq 1$.
Since $\beta_{w}$, $\alpha_{w}$, and $\Gamma_{w}$ are
structurally stable,
we can find a small neighborhood
$U_{w}$ about $c_{w}$ with $\hbox{diam}(U_{w})\leq 1/2^{k}$ such that
the corresponding domains $D_{wn}(c)$ and $B_{wn}(c)$ can be constructed
for $P_{c}$, $c\in U_{w}$. Let
$$
W_{wn} =\{ c\in {\bf C} \; |\; F_{w,c}(0) \in B_{wn}(c)
\}.
$$
The diameter $\hbox{diam}(B_{wn}(c))$ tends to zero and
the set $B_{wn}(c)$ approaches to $\overline{\beta_{w}}(c)$
uniformly on $U_{w}$ as $n$ goes to infinity.
Since the equation $F_{w,c}(0)-\overline{\beta}_{w}(c)=0$
has a unique solution $c_{w}$, the Rouch\'e Theorem implies that
$F_{w,c}(0)-x=0$ for any $x$ in $B_{wn}(c)$ and for $n$ large has a
unique solution which is close to $c_{w}$. Therefore,
there is an integer $N_{w}\geq 0$ such that for $n\geq N_{w}$, $W_{wn}\subseteq
U_{w}$. Thus $\hbox{diam}(W_{wn})\leq 1/2^{k}$.
Because $W_{wn}$ for $n\geq N_{i_{1}}$ is bounded by external rays and
an equipotential curve of ${\cal M}$, 
$$
\tilde{\cal M}_{wn}=W_{wn}\cap {\cal M}$$ 
is connected.

For each $c$ in $W_{wn}$, $n\geq N_{w}$,
let $C_{wn}(c) =F^{-1}_{w,c}(B_{wn}(c))$.
Then
$$
F_{wn,c}=F^{\circ (n+1)}_{w,c} : \CC_{wn}(c) \rightarrow
\DD_{w0}(c)
$$
is a quadratic-like map. Moreover,
$$\{ F_{wn,c}: \CC_{wn}(c) \rightarrow \DD_{w0}(c) \; |\; c\in W_{wn}\}$$ 
is a full family. Thus
$W_{wn}$ contains a copy ${\cal M}_{wn}\subset \tilde{\cal M}_{wn}$ of
the Mandelbrot set ${\cal M}$. For
$$
c\in \Upsilon_{k+1} = \cup_{w} \cup_{i_{k}\geq
N_{w}} {\cal M}_{wi_{k}},
$$
$P_{c}$ is $(k+1)$-{\sl times} renormalizable where
$w=i_{0}i_{1}\ldots i_{k-1}$ runs over all sequences of integers of length
$k$ in the induction.

We have thus constructed a  subset
$\Upsilon (-2) =\cap_{k=1}^{\infty} \Upsilon_{k}$
such that for each $c\in
\Upsilon (-2)$, $P_{c}$ is infinitely renormalizable and such that
$-2$ is a limit point of $\Upsilon (-2)$. For each $c\in \Upsilon
(-2)$, there is a corresponding sequence $w_{\infty} =i_{0}i_{1}\ldots
i_{k}\ldots$ of integers such that
$$
\{c \} =\cap_{k=0}^{\infty}W_{i_{0}\ldots i_{k}}.
$$
Since 
$$
\tilde{\cal M}_{i_{0}\ldots i_{k}}=W_{i_{0}\ldots i_{k}}\cap
{\cal M}
$$
is connected, 
$$
\{ W_{i_{0}\ldots i_{k}}\}_{k=0}^{\infty}
$$ 
is a basis of connected neighborhoods of the Mandelbrot ${\cal M}$ at $c$. In other
words, the Mandelbrot set ${\cal M}$ is locally connected at $c$.

We prove further,
after understood the above construction,
that for any Misiurewicz point, there is a subset $\Upsilon (c_{0})$
such that every $c$ in $\Upsilon (c_{0})$ is infinitely renormalizable, such
that the Mandelbrot set ${\cal M}$ is locally connected at every point
$c$ in $\Upsilon (c_{0})$, and such that $c_{0}$ is a limit point of
$\Upsilon(c_{0})$ as follows.

Suppose $c_{0}$ is a Misiurewicz point in ${\cal M}$.
Then there is an integer $m>1$ such that $p=P_{c_{0}}^{\circ m}(0)$ is a
repelling periodic point of period $k\geq 1$. Let $\alpha$ be the separate
fixed point of $P_{c_{0}}$. Without loss of generality, we assume that
$P_{c_{0}}$ is non-renormalizable. (If $P_{c_{0}}$ is renormalizable,
it must be
finitely renormalizable. We would then take $\alpha$ as the separate
fixed point of the last renormalization of $P_{c_{0}}$ (see \cite{ji})).
Let $\Gamma$ be the union of a cycle of external rays landing at $\alpha$. 
Let $\gamma$ be a fixed equipotential curve of $P_{c_{0}}$. Using $\Gamma$ and
$\gamma$, we can construct the two-dimensional Yoccoz puzzle as follows
(see \cite{ji} for the notation). Let $C_{-1}$ be the domain bounded by $\gamma$. Then $\Gamma$ cuts
$C_{-1}$ into a finite number of closed domains. Let $\eta_{0}$ denote the
set of these domains. Let $\eta_{n} =P_{c_{0}}^{-n} (\eta_{0})$.
Let $C_{n}$ be the member of $\eta_{n}$ containing $0$ for $n\geq
0$.

Let
$$ 
p\in \cdots \subseteq D_{n}(p) \subseteq D_{n-1}(p) \subseteq \cdots
\subseteq D_{1}(p) \subseteq D_{0}(p)
$$
be a $p$-end, where $D_{n}(p) \in \eta_{n}$. Let
$$ 
c_{0}\in \cdots \subseteq E_{n}(c_{0}) \subseteq E_{n-1}(c_{0})
\subseteq \cdots \subseteq E_{1}(c_{0}) \subseteq E_{0}(c_{0})
$$
be a $c_{0}$-end, where $E_{n}(c_{0}) \in \eta_{n}$.
We have $P_{c_{0}}^{\circ (m-1)}(E_{n+m-1}(c_{0})) =D_{n}(p)$.
Since the diameter $\hbox{diam}(D_{n}(p))$ tends to zero
as $n\rightarrow \infty$ and since $p$ is a repelling periodic point,
we can find an integer $l\geq m$ such that $|(P_{c_{0}}^{\circ
k})'(x)|\geq
\lambda >1$ for all $x\in D_{l}(p)$ and such that 
$$
P_{c_{0}}^{\circ (m-1)}: E_{l+m-1}(c_{0}) \rightarrow D_{l}(p)
$$ 
is a homeomorphism.
Let $q\geq 0$ be the integer such that
$$
f=P_{c_{0}}^{\circ q} : D_{l}(p) \rightarrow C_{r_{0}}
$$
is a homeomorphism, where $C_{r_{0}}$ is the domain containing
$0$ in $\eta_{r_{0}}$, $r_{0}\geq 0$.
There is an integer
$r>r_{0}$ such that $r+q>l$ and such that $B_{0}=f^{-q}(C_{r})$ does not
contain $p$. Thus
$$
P_{c_{0}}^{\circ q} : B_{0} \rightarrow C_{r}
$$
is a homeomorphism. Define
$$
B_{n}= \Big( P_{c_{0}}^{\circ nk}|D_{l+nk}(p)\Big)^{-1}(B_{0})
\subseteq D_{l+nk}(p)
$$
for $n\geq 1$. Note that $B_{n}$ is in $\eta_{r+q+nk}$. Then
$$
P_{c_{0}}^{\circ (q+nk)} : B_{n} \rightarrow C_{r}
$$
is a homeomorphism.

Since $\alpha$, $p$, $\Gamma$, $C_{r}$, $D_{n}$, and $B_{n}$, for $n\geq 0$, are
structurally stable, they can be constructed for $P_{c}$ as long as $c$
near $c_{0}$. Let $U_{0}$ be a neighborhood about $c_{0}$ with
$\hbox{diam}(U_{0})\leq 1$ such that the corresponding points
$\alpha(c)$ and $p(c)$ and the corresponding sets $\Gamma (c)$,
$C_{r}(c)$, $D_{n}(c)$,
and $B_{n}(c)$, for $n\geq 0$, are all preserved for $c\in U_{0}$.
Moreover, as $n$ goes to infinity, the diameter $\hbox{diam}(B_{n}(c))$
tends to zero and the set $B_{n}(c)$ approaches to $p(c)$ uniformly on $U_{0}$.
Let
$$
W_{n}=W_{n}(c_{0})=\{ c \in {\bf C} \; |\; P_{c}^{m} (0) \in
B_{n}(c)\}.
$$
Since the equation $P^{\circ m}_{c}(0)-p(c)=0$ has a unique solution
$c_{0}$ in $U_{0}$ (see \cite{dh2}), the Rouch\'e
Theorem implies that the equation
$P^{\circ m}_{c}(0)-x=0$ has a unique solution for any $x$ in $B_{n}(c)$
and large $n$ and that the solution is close to $c_{0}$.
Therefore, there is an integer
$N_{0}=N_{0}(c_{0})>0$ such that $W_{n}\subset U_{0}$. Thus
$\hbox{diam}(W_{n})\leq 1$ for all $n\geq N_{0}$.
Since $W_{n}$ is bounded by external rays and an equipotential
curve of the Mandelbrot set ${\cal M}$, 
then 
$$
\tilde{\cal M}_{n} ={\cal M}\cap W_{n}
$$ 
is connected (see Fig. 5).

\begin{figure}
\centerline{\psfig{figure=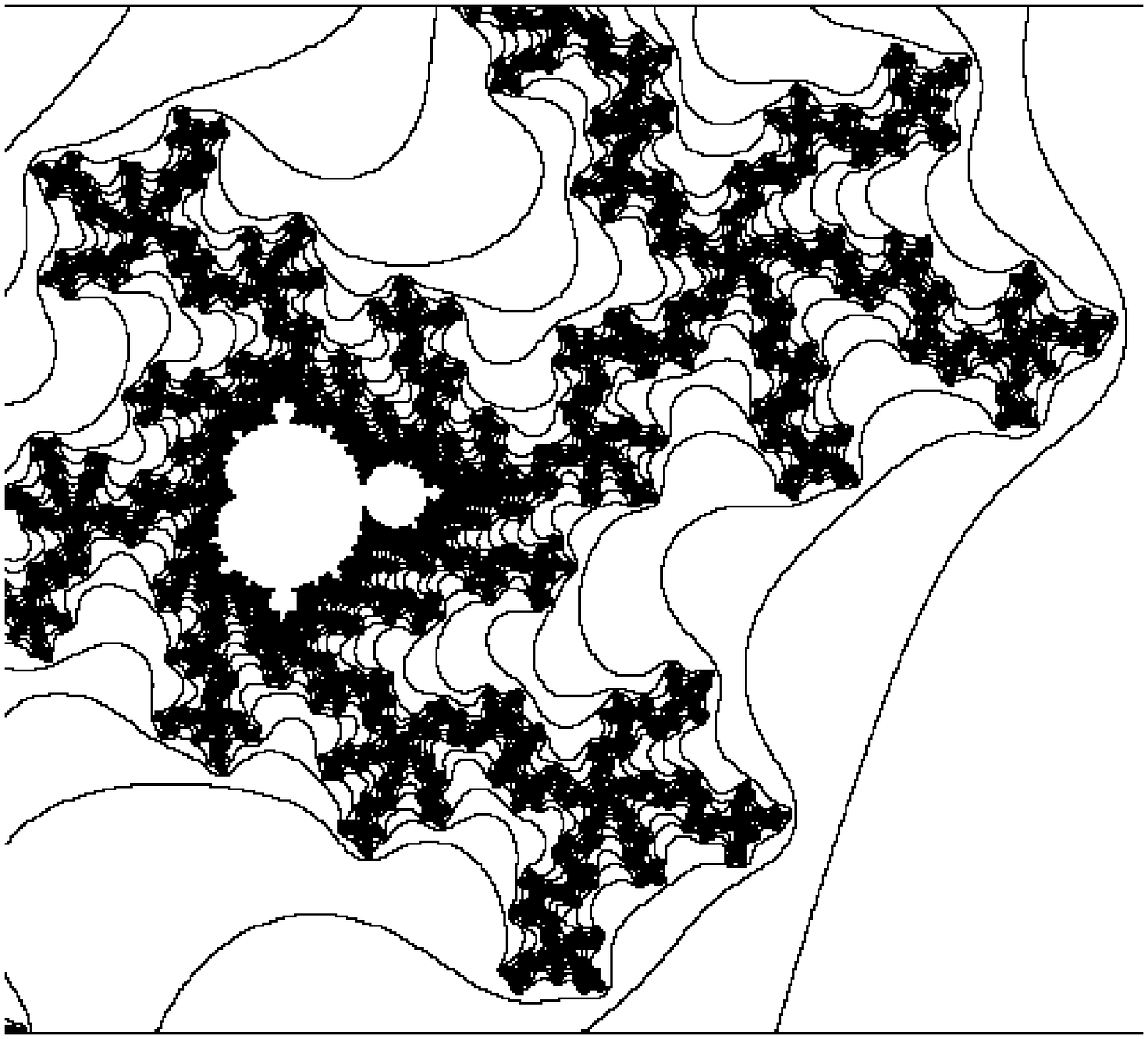,height=3in}}
\caption{A small copy of the Mandelbrot set and hairs around it}
\end{figure}

\vskip5pt
For any $c\in W_{n}$, $n\geq N_{0}$, let $R_{n}(c)$ be the pre-image of
$B_{n}(c)$ under the map
$$P_{c_{0}}^{\circ (m-1)}: E_{l+m-1}(c)\rightarrow D_{l}(p,c)$$
and let $C_{m+r+q+nk}(c)=P_{c}^{-1}(R_{n}(c))$. Then $C_{m+r+q+nk}(c)$
is the domain containing $0$ in $\eta_{m+r+q+nk}$ and
$$F_{n,c}=P_{c}^{\circ (q+nk+m)}:
\CC_{m+r+q+nk}(c)\rightarrow \CC_{r}(c)$$
is a quadratic-like map. Moreover,
$$ \{ F_{n,c}: \CC_{m+r+q+nk}(c)\rightarrow \CC_{r}(c) \; |\; c\in W_{n}\}$$
is a full family. Thus $W_{n}$ contains a copy ${\cal
M}_{n}={\cal M}_{n}(c_{0})$ of the Mandelbrot set ${\cal M}$ where
${\cal M}_{n}\subset
\tilde{\cal M}_{n}$. For any $c\in \cup_{n\geq N_{0}} {\cal M}_{n}$, $P_{c}$ is
once renormalizable.
Now repeat the induction we did for $-2$. We construct a
subset
$\Upsilon (c_{0})$ such that every $c\in \Upsilon (c_{0})$ is
infinitely renormalizable, such that ${\cal M}$ is locally connected at
every $c\in \Upsilon (c_{0})$, and such that $c_{0}$ is a limit point
of $\Upsilon (c_{0})$.

Let $\Upsilon =\cup_{c_{0}}\Upsilon (c_{0})$ where $c_{0}$ runs over all
Misiurewicz points in ${\cal M}$. Then every $c\in \Upsilon$
is infinitely renormalizable and ${\cal M}$ is locally connected at
every $c\in \Upsilon$. Since the set of Misiurewicz points is 
dense in $\partial {\cal M}$, the set $\Upsilon$ is 
dense in $\partial {\cal M}$. This completes the proof of the theorem.

\vskip5pt
{\bf Remark 1.} Eckmann and Epstein \cite{ee} and Douady and
Hubbard \cite{dh3} estimated the size of ${\cal M}_{n}$ in the construction.
Since $\tilde{\cal M}_{n}\setminus {\cal M}_{n}$ contains hairs (see Fig. 5)
and may destroy the local connectivity of ${\cal M}$, we must estimate the
size of $\tilde{\cal M}_{n}$.

\vskip5pt
{\bf Remark 2.} 
For $c$ in $\Upsilon$, one can prove that 
the Julia set $J_{c}$ is locally connected at the critical point $0$
as well as at the grand critical orbit of $P_{c}$. 
But since we did not show that $P_{c}$ is unbranched 
we could not conclude that $J_{c}$ is locally connected from the
result in \cite{ji}. However, we constructed a similar subset  
$\tilde{\Upsilon}$ of the Mandelbrot set ${\cal M}$ in \cite{ji} 
which is also dense on the boundary of the Mandelbrot set ${\cal M}$ 
such that $P_{c}$ is unbranched infinitely renormalizable and has 
the complex a priori bounds for every point $c$ of this subset. 
Thus the Julia set $J_{c}$ of $P_{c}$ for $c$ in $\tilde{\Upsilon}$ 
is locally connected according to the result in \cite{ji}.

\end{document}